\documentclass[12pt]{article}
\usepackage{amssymb}

\title{Decreasing families of dynamically determined intervals in the power-law family}
\author{Waldemar Pa\l uba\thanks{Partially
                         supported by a KBN grant no.
                         2 PO3A 010 22. }\\
                         \small Institute of Mathematics\\
                         \small Warsaw University\\
                         \small Banacha 2\\
                         \small 02-097 Warsaw, Poland\\
                         \small e-mail: paluba@mimuw.edu.pl
                         }
\date{}

\makeatletter \@addtoreset{equation}{section}

\makeatother

\newtheorem{prop}{Proposition}[section]
\newtheorem{thm}{Theorem}[section]
\newtheorem{lemma}{Lemma}[section]
\newtheorem{cor}{Corrolary}[section]
\begin{document}

\maketitle

\begin{abstract}
We study the rate of growth of ratios of intervals  delimited by
the post-critical orbit of a map in the quasi-quadratic family
$x\mapsto -|x|^\alpha +a.$ The critical order $\alpha$ is an
arbitrary real number $\alpha>1.$ The range of the parameter $a$
is confined to an interval $(1,a_{\alpha})$ of length depending on
the critical order. We prove that in every power-law family there
is a unique parameter $p_{\alpha}$ corresponding to the kneading
sequence $RLRRRLRC.$ Subsequently, we obtain  monotonicity results
concerning ratios of  all intervals labeled by infinite
post-critical orbit in the case of the kneading sequence $RLRL...$
This extends the results from \cite{P}, via refinement of the
tools based on special properties of power-law mappings in
non-euclidean metric.
\end{abstract}


\section{Introduction}\label{intr}
 In this paper we continue our work done in \cite{P} on families
of
 unimodal quasi-quadratic maps of the form $f_a(x)= -|x|^{\alpha} + a,$
  with a real parameter $a$ and an arbitrary -- in general
  non-integer -- fixed exponent $\alpha >1\, .$

 The problem of  monotone behaviour of the dynamics in such a
family has  been first successfully solved for the strictly
quadratic case $\alpha =2\, .$ The tools initially  developed for
the quadratic case (see e.g  \cite{GS2}, \cite{Ly}, \cite{Sh},
also \cite{Ts}, and an independent attempt, partly relying on
numerical evidence in \cite{DJS}) were broadly generalized in the
work of Kozlovski-Shen-van Strien, see \cite{KSS}. There are also
very interesting recent results by G.Levin,  concerning uniqueness
of appearance of periodic orbits of given multiplier in the
quadratic family $z^2+c$. Not only was he able to give a simple
proof of Douady-Hubbard-Sullivan theorem (cf.\cite{L1}), but he
could continue somewhat beyond the hyperbolic domains in the
Mandelbrot set also, see \cite{L2}. In this work, we focus on
questions closely  related to  these of Levin's, though only
orbits of periods $2$ or $4$ appear here. In return, working with
real variable tools, we can do all critical degrees, integer or
non-integer, indiscriminately.

 Despite of  a
great deal of progress achieved in the aforementioned papers, and
in other works as well, virtually all those developments are
inherently limited to the case of integer critical degrees.
Non-integers clearly require a fresh and different approach.
 For any real number $\alpha >1\,$ the power-law map $\, x\mapsto |x|^{\alpha}\,
$ has negative Schwarzian derivative, and hence it expands the
non-euclidean lengths. This observation has long become one of the
key tools in one-dimensional dynamics. However, the power-law is
not just {\it a} negative Schwarzian map. It is a homogeneous map,
and in the Poincar\'e metric with the element
 $\,\frac{dt}{t}\,$ on the positive half-line $\,(0,\infty),$ it
is nothing but a linear map acting as multiplication by the
coefficient
 $\alpha$, once we set the origin of the Poincar\'e coordinates at $1.\,$
This simple fact is rather hard to make use of in a direct way,
but carries some strong consequences  that can be applied in a
dynamical setting.

 In our previous paper on this subject (see \cite{P}), we
introduced the technique of indirect use of linearity of the
power-law map in the non-euclidean metric and exemplified its
usefulness  in dynamics. There, we studied maps in the
one-parameter quasi-quadratic family $f_a\,$ with the kneading
sequence $RRR\ldots\,$, that is for the value of the parameter
$a\,$ smaller than $1\, .$ For the infinite decreasing family of
intervals with endpoints labeled by the successive points of the
post-critical orbit we proved that the ratio of any two such
intervals is a function monotone in parameter $a\, .$ This means,
we studied the situation which arises before the orbit of the
critical point becomes a super-stable orbit of period $2\, .$ It
is clear that this period $2$ super-sink situation arises only
once in our family.

 In the current work, we further develop our tools in order
to examine the case of some parameters greater than $1$, where the
length of the interval $(1,a_{\alpha})$ to which those parameters
are confined depends on the critical order $\alpha
>1\, .$ In particular we shall be able to deal with kneading
sequences of the form $RLRL\ldots\,$,  proving monotonicity of the
ratios, in the respective decreasing families, of intervals
delimited by the post-critical orbit.  As a step in the build-up
of the above techniques, we shall also establish uniqueness of the
period $8$ super-sink, corresponding to the kneading sequence
$RLRRRLRC$, in every power-law family (even when it does not admit
a holomorphic extension!); uniqueness of the period $4$ super-sink
$RLRC$ is elementary and follows along the way.

 \section{Notation and preliminaries}\label{prelim}
 To begin with, we set some notation in conformance with that of
\cite{P}. The names {\it non-euclidean} and {\it Poincar\'e} we be
used interchangeably.

 The Poincar\'e coordinate of a point $x$ in an  oriented, open
interval $(p,q)$ will be denoted by $$p_{p,q}(x)=\ln
{\frac{x-p}{q-x}}\,
 ,$$ and respectively $p_{p,\infty}(x)=\ln (x-p)$;
  also $p_{p,-\infty}(x)=\ln |x-p| $.

 To single out the non-euclidean metric on the half-line, which
 turns the mapping $h$ into a linear map, we will coin the term
 {\it nonlinearity of an interval} for the length of this interval
 measured in the Poincar\'e metric on $(0,\infty)$. Under this
 convention, the integral of nonlinearity of $h$ over an interval $(p,q)$
 equals, up to a multiplicative constant, to the {\it
 nonlinearity} of the domain of integration.

 Given an orientation preserving  homeomorphism
  $\varphi : (p,q)\to (r,s)\,$ we shall observe the `bar'
  notation for its counterpart in the
non-euclidean coordinates, i.e. the mapping
$\overline\varphi_{p,q} : \mathbb{R}\to\mathbb{R}\,$ defined by
the formula
$$\overline\varphi_{p,q}(t)=p_{r,s}
 \left(\varphi\left(p^{-1}_{p,q}(t)\right)\right)\, .$$
The {\it non-euclidean push }  of $\varphi$ at a point $x\in
(p,q)$ is, by definition, the quantity
 $$p_{r,s}\left( \varphi(x) \right)-p_{p,q}(x)\, .$$
By the {\it strength} of a push we mean its absolute value.

By $\overline\varphi^+_{p,q}$ and $\overline\varphi^-_{p,q}$ we
denote the finite or infinite limits
$$\overline\varphi^+_{p,q} = \lim_{t\rightarrow +\infty}\left( \overline\varphi(t)-t \right)$$
$$\overline\varphi^-_{p,q} = \lim_{t\rightarrow -\infty}\left( \overline\varphi(t)-t \right),$$
provided they exist.

 When $\varphi$ is the restriction of the homogeneous map
  $h (x)=x^\alpha$ to an interval $(p,q)\subseteq (0,\infty)$
 we shall always put $h$ or $\overline h$ in place of $\varphi$ or
 $\overline\varphi$ respectively.

  For a fixed exponent $\alpha >1,$ let $f_a=-|x|^{\alpha}+a\,$,
and the successive points of the orbit of the critical point will
simply be denoted by $n_a=f_{a}^n(0)\, .$ Moreover, homogeneity of
the power-law map allows for the linear change of coordinates, $\,
n_a\mapsto n_a/1_a\, $, so that we can set $0_a=0,\, 1_a=1\,$ and
the dependence on the parameter $a$ turns into the dependence on
the value of $2_a$ in these new coordinates. For $a>1$ this
rescaled value of $2_a$ is  in the interval $(-1,0)$  -- so long
as the post-critical orbit does not escape to infinity --  and the
quantity $p_{0,1}(|2_a|)$, which for obvious reason we will denote
by $\bar a\,$,  is increasing simultaneously with $a $. Throughout
this work, this very quantity will be chosen as our new parameter,
and it is always tacitly assumed that the rescaling $n_a\mapsto
n_a/1_a\,$ has been done.

We now record several  observations concerning one-dimensional
non-euclidean coordinates. Below, they are  stated as
propositions,   verifiable by  elementary calculations derived
directly from the definition of the Poincar\'e metric.
\begin{prop}\label{prop2.1}
For any $x\in (-\infty ,0)$ the following two Poincar\'e
coordinates coincide
$$p_{0,-\infty}(x)= p_{x,1}(0) .$$
\end{prop}
{\it Proof}. We have
$p_{0,-\infty}(x)=\ln{(-x)}=\ln{\frac{0-x}{1-0}}=p_{x,1}(0)$.
 \hfill$\Box$

\begin{prop}\label{prop2.2}
For any $x\in (-1,0)$ the following two Poincar\'e coordinates
coincide
$$p_{1,-1}(x)=p_{x,-x}\left(p_{x,1}^{-1}\left(p_{0,-1}(x) \right) \right)
.$$
\end{prop}
{\it Proof}. The identity in question is tantamount to
$\ln{\frac{x-1}{-1-x}}=\ln{\frac{c-x}{-x-c}}$, or
\begin{equation}\label{eq2.1}
\frac{x-1}{1+x}=\frac{c-x}{x+c} ,\end{equation}
 where $c=p_{x,1}^{-1}(p_{0,-1}(x)) $, i.e. $p_{x,1}(c)=p_{0,-1}(x) .$
But this last equality means $\frac{c-x}{1-c}=\frac{x}{-1-x}$, and
further $c=x^2$, so that (\ref{eq2.1}) follows.
  \hfill  $\Box$

Given a point $x\in (-\infty,0), $ we then pick a   point $y\in
(x,0) .$ We shall  let the point $x$ vary, by which we mean a
choice of another point $\tilde x \in (-\infty ,0).$ The
discrepancy in the non-euclidean coordinate will be denoted by
$$\Delta \vartheta = p_{0,-\infty}(\tilde x) - p_{0,-\infty}(x).$$
A broader version of Proposition \ref{prop2.1} is the following
\begin{prop}\label{prop2.3}
In the above notation  we have
$$  p_{\tilde x, 0}\left(p_{\tilde x,1}^{-1}
 \left(p_{x,1}(y)+\Delta \vartheta \right) \right) -p_{x,0}(y) =
  p_{1,-\infty}(\tilde x)-p_{1,-\infty}(x)  .$$
\end{prop}
{\it Proof}. We have $\Delta\vartheta=\ln{ \frac{\tilde x}{x}}\, $
and $\, p_{1,-\infty}(\tilde x)-p_{1,-\infty}(x)=\ln{\frac{\tilde
x -1}{x-1}}\,$. Denote  $ c=p^{-1}_{\tilde
x,1}(p_{x,1}(y)+\Delta\vartheta)$, a point characterized by
\begin{equation}\label{eq2.2}
\frac{c-\tilde x}{1-c}=\frac{y-x}{1-y}\cdot\frac{\tilde x}{x}\, .
\end{equation}
We will be done once we show $\frac{c-\tilde
x}{-c}\cdot\frac{-y}{y-x}=\frac{1-\tilde x}{1-x}\, $, or
\begin{equation}\label{eq2.3}
\left( 1+\frac{\tilde x}{c-\tilde x} \right)\left(\frac{y-x}{y}
\right)=\frac{1-x}{1-\tilde x}\, .
\end{equation}
From (\ref{eq2.2}) we get $\frac{\tilde x}{c-\tilde
x}=\frac{1-y}{y-x}\cdot\frac{x}{1-\tilde x}+\frac{\tilde
x}{1-\tilde x}\, $ and (\ref{eq2.3}) can now be checked
immediately.
   \hfill$\Box$\\

Proposition \ref{prop2.1} is what we get of Proposition
\ref{prop2.3}, in place of subtracting two infinite terms,  when
we set $y=0$. We generalize Proposition \ref{prop2.2} in a similar
way. Suppose we are given a point $x\in (-1,0),$ and a point $y\in
(x,-x)$. Again, we let the point $x$ vary by choosing a new point
$\tilde x\in (-1,0)$. The discrepancies in the appropriate
Poincar\'e coordinates of the two points will be denoted by
$$\Delta t=p_{1,-1}(\tilde x)-p_{1,-1}(x)\, , $$
and by
$$\Delta\theta =p_{0,-1}(\tilde x)-p_{0,-1}(x) $$ respectively.
A statement parallel to Proposition \ref{prop2.3} is the following
\begin{prop}\label{prop2.4}
In the above notation we have
 $$p_{\tilde x,-\tilde x}\left(p_{\tilde
 x,1}^{-1}\left(p_{x,1}(y)+\Delta\theta
  \right) \right) - p_{x,-x}(y)=\Delta t .$$
   \end{prop}
{\it Proof}. The point $c=p^{-1}_{\tilde
x,1}(p_{x,1}(y)+\Delta\theta)\,$ satisfies $\, \frac{c-\tilde
x}{1-c}=\frac{y-x}{1-y}\cdot \frac{\tilde x}{1+\tilde
x}\cdot\frac{1+x}{x}\, $, which can be transformed into
\begin{equation}\label{eq2.4}
\frac{\tilde x}{c-\tilde x}=\frac{1-y}{y-x}\cdot\frac{1+\tilde
x}{1-\tilde x}\cdot\frac{x}{1-x}+\frac{\tilde x}{1-\tilde x}\, .
\end{equation}
We will be done if we show that $\frac{c-\tilde x}{\tilde
x+c}\cdot\frac{x+y}{y-x}=\frac{\tilde x -1}{1+\tilde
x}\cdot\frac{1+x}{x-1}\, $, which is the same as
 \begin{equation}\label{eq2.5}
 \frac{\tilde x +c}{c-\tilde x}=\frac{y+x}{y-x}\cdot\frac{1+\tilde x}{\tilde x
 -1}\cdot\frac{x-1}{1+x}\, .
 \end{equation}
 Since $\frac{\tilde x +c}{ c-\tilde x}=1+\frac{2\tilde x}{c-\tilde
 x}\,$, equation (\ref{eq2.5}) follows immediately from (\ref{eq2.4}).
   \hfill$\Box$

\begin{prop}\label{prop2.5}
Suppose $x,x'\in(0,1)$ and $y\in(x,1).$ Let $y'$ be such a point
in $(x',1)$ that $p_{x',1}(y')=p_{x,1}(y).$ Then
\begin{equation}\label{kkk}
p_{1,0}(x')-p_{1,0}(x)=p_{y',0}(x')-p_{y,0}(x).
\end{equation}
\end{prop}
{\it Proof.} The point $y'$ is chosen in such a way that
$\frac{y'-x'}{1-y'}=\frac{y-x}{1-y},$ or
$\frac{1-x'}{1-y'}=\frac{1-x}{1-y}.$ Identity (\ref{kkk}) is now
immediate.\hfill$\Box$

\section{The period $4$ super-sink}\label{period4}
In this short section we describe the behavior of the point $4_a$
when we let the parameter $\bar a$ vary in such a range, that
$3_a\in (0,1)$ and the point $4_a$ stays within the interval
$(2_a,-2_a) .$

Let a positive number $t$ be the Poincar\'e coordinate of $2_a$ in
the oriented interval $(1,-1)$, and we set
$$g(t)=p_{2_a,-2_a}(4_a) .$$
The following theorem holds true.
\begin{thm}\label{thm3.1}
The inverse function $g^{-1} : \mathbb{R}\to\mathbb{R}_+$ is
strictly increasing, and $g'(t)>1$. In particular, the value
$g(t)=0$, corresponding to  the  super-stable orbit with the
kneading sequence RLRC is assumed only once.
\end{thm}
 {\it Proof}. Consider a pair of admissible  parameter values $\bar
 a$ and $\bar a'$, i.e. such that the orbits $2_a,\, 3_a,\, 4_a\,$
(and respectively $2_{a'},\, 3_{a'},\, 4_{a'}$) satisfy the
restrains on the dynamics we set above. Then
\begin{equation}\label{deltat}
\Delta t=p_{1,-1}(2_{a'})-p_{1,-1}(2_a),\mbox{ while } \Delta g=
p_{2_{a'},-2_{a'}}(4_{a'})-p_{2_a,-2_a}(4_a). \end{equation}
 Applying Proposition \ref{prop2.4} to this case we get
 \begin{equation}\label{pierwsze}
p_{2_{a'},-2_{a'}}(p^{-1}_{2_{a'},1}(p_{2_a,1}(4_a)+(p_{0,-1}(2_{a'})-p_{0,-1}(2_a))))-
p_{2_a,-2_a}(4_a)=\Delta t,
\end{equation}
so, because of monotonicity of the coordinate functions, we only
need to establish that
 \begin{equation}\label{nabla}
p_{2_{a'},1}(4_{a'})-p_{2_a,1}(4_a)>p_{0,-1}(2_{a'})-p_{0,-1}(2_a).
\end{equation}
This inequality becomes clear once we split the procedure leading
from point $2_a$ (respectively $2_{a'}$) to $4_a$ (respectively to
$4_{a'}$) into three steps. In the first step, we act on the
interval $ (0,-1)$ by the restriction of the power-law map. Thus,
due to negative Schwarzian derivative, the initial discrepancy
($p_{0,-1}(2_{a'})-p_{0,-1}(2_a)$) in the Poicar\'e coordinates
gets increased. So we see that
$$p_{1,2_{a'}}(3_{a'})
-p_{1,2_a}(3_a)\geq p_{0,-1}(2_{a'})-p_{0,-1}(2_a).$$

 In the
second step, we turn the interval $(0,-1)$ over, onto the interval
$(1,2_a)$, or  onto $(1,2_{a'})$ respectively, and then we
truncate the image at the critical point 0. This cut-off only
increases the Poincar\'e coordinate of every point, which after
the turnover landed in $(1,0)$,  because we now read the
Poincar\'e coordinate in the  interval $(1,0)$ rather than in a
larger domain $(1,2_a)$, or $(1,2_{a'})$ respectively.
 Moreover, the increase in the Poincar\'e coordinate  inflicted
by cutting the domain interval short, is  in the case of point
$3_{a} $ smaller then in the case of $3_{a'} $.
 This is so, because
  the endpoint $2_{a}$  is closer to the critical point,
  while the endpoint $ 2_{a'}$ is further away to the left,
     so of two corresponding points with identical  Poincar\'e coordinate
      within the
respective domain intervals (with the other endpoint at $1$), the
gain  in the latter situation is larger than in the former. But
instead of equal coordinates, we have even better inequality
$p_{1,2_{a'}}(3_{a'})>p_{1,2_a}(3_a) $, which  further enlarges
the gain.
 Thus, in this second step,
 made of the turnover
 followed by truncation,  the initial discrepancy grows even further and so
\[p_{1,0}(3_{a'})-p_{1,0}(3_a)>p_{1,2_{a'}}(3_{a'})
-p_{1,2_a}(3_a).\]

In the last step, we again act by a negative Schwarzian map
stretching the discrepancy between the Poincar\'e coordinates yet
further, and finally we make the turnover onto $(1,2_a)$, and
respectively onto $(1,2_{a'})$, to arrive at (\ref{nabla}).
Therefore $\Delta g > \Delta t $ and the proof is complete.
\hfill$\Box$
\section{The period $8$ super-sink}\label{period8}
In the previous section we have established that, when we vary the
parameter $\bar a\,$, the position of the point $4_a$ within the
interval $(2_a, -2_a)$ changes monotonically,  with the derivative
greater than $1$. It   clearly follows from the proof,  that this
derivative actually stays bounded away from $1$, in a way that
depends on the critical order $\alpha\,$. In
  section \ref{rlrl} we will study in detail the case  of $p_{2_a,-2_a}(4_a)<0$, and
  describe the behavior of the intervals delimited by  the post-critical orbit with
  the kneading sequence $RLRL\ldots\, $.

In here, we will focus on these admissible parameters $\bar a\,$,
for which \linebreak $p_{2_a,-2_a}(4_a)>0\,$ and
$p_{4_a,-4_a}(8_a)\leq 0\, $,
 i.e. we are past the (unique) parameter corresponding to $RLRC$,
  but we do not cover the critical point yet another time. From now on, we are
  making our choice of the parameter subject to this restriction.
 We shall see that, as long as the above
 condition on the dynamics is satisfied, the movement of the point
  $8_a$ is also monotone in
parameter, and in the non-euclidean metric in $(4_a, -4_a)$ this
point moves with the derivative strictly positive. It will follow
that the $RLRRRLRC$ super-stable orbit appears uniquely in every
power-law family. It is a subject of an ongoing work, that goes
beyond the scope of this paper, to  examine whether a claim
analogous to that of Theorem \ref{thm3.1} can be fully extended to
larger set of parameters.

In our current case,  the scheme of the argument we used to prove
Theorem \ref{thm3.1}, alone will  not suffice, and a more delicate
technique must be employed. Yet,
  some understanding of the way Poincar\'e coordinates vary
  remains an important component.
Since, due to the more intricate dynamics, the required property
of the non-euclidean coordinates becomes less self-evident, we
state it as a separate lemma. The points $\,x,y,z\, $ below will
correspond to the points $\, 2_a,4_a,8_a\,$ of the post-critical
orbit. The origin of the summands, which do not have equivalent in
the statement of Proposition \ref{prop2.4} will be explained
later, in the course of the proof of Theorem \ref{thm4.1} below.
Here, we only indicate that the last term has to do with the limit
strength of a non-euclidean push.
\begin{lemma}\label{osemka}
Suppose we are given two triples of points, $\,(x,y,z)$ and
$(\tilde x,\tilde y,\tilde z )$, satisfying the following
conditions:
\begin{description}
\item[(i)] $ x,\tilde x\in (0,-1)\, $  and $ \, p_{0,-1}(\tilde x)>
p_{0,-1}(x)\,$,
\item[(ii)] $ y\in (0,-x)$, $\, \tilde y\in (0,-\tilde x)\, $ and
$\, p_{\tilde x,1}(\tilde y)\geq p_{x,1}(y)+\left(p_{0,-1}(\tilde
x )-p_{0,-1}(x) \right)\,$,
\item[(iii)] $z\in (y,0] $, $\, \tilde z\in (\tilde y,-\tilde y)\,$ and
$\, p_{\tilde y,\tilde x}(\tilde z)\geq p_{y,x}(z)+(p_{0,-\tilde x
}(\tilde y)-p_{0,-x}(y) )+\ln{\frac{\tilde y -\tilde x}{y-x}}-
  \left((p_{\tilde x,1}(\tilde y)-p_{x,1}(y) ) -
   (p_{0,-1}(\tilde x)- p_{0,-1}(x) )  \right) \, $.
\end{description}
Then $\, p_{\tilde y,-\tilde y}(\tilde z)>p_{y,-y}(z)\,$.
\end{lemma}
{\it Proof}. It is immediate to check that for arbitrary $y,\tilde
y\in (0,1)$ one has
\begin{equation}\label{eq4.9}
p_{\tilde y,-1}(0)=p_{y,-1}(0)+\left(p_{0,1}(\tilde y)-p_{0,1}(y)
\right)+\ln{\frac{1+\tilde y}{1+y} } - \left(p_{-1,1}(\tilde
y)-p_{-1,1}(y) \right).
\end{equation}
 We now assume $\tilde y>y$, and allowing $z\neq 0$ we verify,
   that for any $z\in [0,y)$
 the following generalization of (\ref{eq4.9}) holds
\begin{eqnarray}\label{eq4.10}
p_{\tilde y,-\tilde y}(p^{-1}_{\tilde y,-1}(p_{y,-1}(z)+
(p_{0,1}(\tilde y)-p_{0,1}(y) )+\ln{\frac{1+\tilde y}{1+y}
}\nonumber\\
 - (p_{-1,1}(\tilde y)-p_{-1,1}(y) ) ) ) \geq p_{y,-y}(z) .
\end{eqnarray}
In order to see this, notice that $$(p_{0,1}(\tilde
y)-p_{0,1}(y))+\ln\frac{1+\tilde y}{1+y}-(p_{-1,1}(\tilde
y)-p_{-1,1}(y) )=\ln\frac{\tilde y}{y}\, ,  $$ and denote
$c=p^{-1}_{\tilde y,-1}(p_{y,-1}(z)+\ln\frac{\tilde y}{y}) $,
which means $\frac{c-\tilde y}{-\tilde
y-c}=\frac{z-y}{1-z}\cdot\frac{\tilde y}{y} $, or
\begin{equation}\label{eq4.11}
\frac{\tilde y}{c-\tilde y}=\frac{y}{1+\tilde
y}\cdot\frac{1+z}{z-y}-\frac{\tilde y}{1+\tilde y} .
\end{equation}
 We will be done if we show that $\frac{c-\tilde y}{-\tilde y-c}\geq
  \frac{z-y}{y-z} $, being equivalent to $\frac{2\tilde y}{c-\tilde y}+1\geq
  \frac{y+z}{z-y}  $, or $\frac{\tilde y}{c-\tilde y}\geq\frac{y}{z-y}
  $. The last inequality follows from (\ref{eq4.11}), once we recall $\tilde y\geq
  y$.

  In the next step we  extend formula (\ref{eq4.10}), allowing $\tilde x\neq
  -1$.  Assuming $1> -\tilde x>\tilde y>y>z\geq 0$, we will now show that
\begin{eqnarray}\label{eq4.12}
p_{\tilde y,-\tilde y}(p^{-1}_{\tilde y,\tilde x}(p_{y,\tilde
x}(z) + (p_{0,-\tilde x}(\tilde y)-p_{0,-\tilde x}(y) )\nonumber
\\
+ \ln\frac{\tilde y -\tilde x}{y-\tilde x}- (p_{\tilde x,1}(\tilde
y)-p_{\tilde x,1}(y) ) ) )>p_{y,-y}(z)  .
\end{eqnarray}
 We emphasize that the inequality in formula (\ref{eq4.12}) is
 always sharp, even for $z=0$.

 This time, we set
\begin{equation}\label{pointc}
 c=p_{y,\tilde
x}(z) + (p_{0,-\tilde x}(\tilde y)-p_{0,-\tilde x}(y) )
+\ln\frac{\tilde y -\tilde x}{y-\tilde x}- (p_{\tilde x,1}(\tilde
y)-p_{\tilde x,1}(y) ),
\end{equation}
 which means
\[ \frac{c-\tilde y}{\tilde x-c}=\frac{z-y}{\tilde x-z}\cdot\frac{\tilde y}{-\tilde x-\tilde y}
\cdot\frac{-\tilde x-y}{y}\cdot\frac{1-\tilde y}{1-y}.\] We
transform this identity into
\[
\frac{\tilde x-\tilde y}{c-\tilde y}-1 =\frac{\tilde
x-z}{z-y}\cdot \frac{\tilde x+\tilde y}{\tilde y}\cdot
\frac{y}{\tilde x +y}\cdot \frac{1-y}{1-\tilde y}
\]
and further into
\[
\frac{\tilde y}{c-\tilde y}=\frac{\tilde y}{\tilde x-\tilde
y}\left[ \frac{\tilde x-z}{z-y}\cdot \frac{\tilde x+\tilde
y}{\tilde y}\cdot \frac{y}{\tilde x +y}\cdot \frac{1-y}{1-\tilde
y}+1\right].
\]
We will be done if we show $p_{\tilde y,-\tilde y}(c)>
p_{y,-y}(z),$ i.e. $\frac{c-\tilde y}{-\tilde y -c}>
\frac{z-y}{-y-z},$ which is equivalent to $\frac{\tilde
y}{c-\tilde y}>\frac{y}{z-y},$ and so it is enough to verify that
\[
\frac{\tilde y}{\tilde x-\tilde y}\left[\frac{\tilde
x-z}{z-y}\cdot \frac{\tilde x+\tilde y}{\tilde y}\cdot
\frac{y}{\tilde x +y}\cdot \frac{1-y}{1-\tilde y}+1\right]
>\frac{y}{z-y}.
\]
This inequality can be rewritten as
\[
\frac{\tilde x -z}{z-y}\cdot \frac{\tilde x+\tilde y}{\tilde x
-\tilde y}\cdot \frac{y}{\tilde x +y}\cdot\frac{1-y}{1-\tilde
y}>\frac{y}{z-y}-\frac{\tilde y}{\tilde x-\tilde y},
\]
or (recall that $y<z,$ $\tilde x <0,$ $\tilde y >0$)
\[
(\tilde x-z)(\tilde x+\tilde y)y(1-y)<(\tilde x +y)(1-\tilde
y)(y\tilde x-z\tilde y),
\]
and further
\[
\tilde x y (1+\tilde x)(\tilde y-y) <z (\tilde y -y)(\tilde
x\tilde y +\tilde x y-\tilde x +y\tilde y).
\]
To conclude, we cancel out $(\tilde y-y),$ and observe that
\begin{equation}\label{www}
\tilde x y +\tilde x\tilde y -\tilde x + y\tilde y > \tilde x
(1+\tilde x).
 \end{equation}
This is so because (\ref{www}) boils down to the inequality \(
\tilde x^2+\tilde x(2-y-\tilde y)-y\tilde y <0,\) which is
elementarily true for all $\tilde x\in(-1,0)$ and $y,\tilde y\in
(0,1).$ For completion of the proof we now consider   an arbitrary
point $x\in (0,\tilde x)$, such that $y<-x$. We consider the
movement of $x$-variable from position $x$ to $\tilde x$ and apply
Proposition \ref{prop2.4} twice, first to the induced movement of
$y$-variable, then to the consequent movement of $z$-variable. By
virtue of that Proposition, we see that points $\hat y$ and $\hat
z$, determined by the identities
$$p_{\tilde x,1}(\hat y)=p_{x,1}(y)+(p_{1,-1}(\tilde x)-p_{1,-1}(x) ) $$
$$p_{\hat y,\tilde x}(\hat z)=p_{y,x}(z)+ (p_{1,-1}(\tilde x)-p_{1,-1}(x) ) $$
satisfy $\hat y<\tilde y$ and $p_{\hat y,-\hat y}(\hat
z)=p_{y,-y}(z)+(p_{1,-1}(\tilde x)-p_{1,-1}(x) ) >p_{y,-y}(z)$.
Thus obviously $\ln\frac{\tilde y-\tilde x}{y-x}>\ln\frac{\tilde
y-\tilde x }{\hat y-\tilde x}$.

 If $\hat z\geq 0$, i.e. $p_{\hat y,-\hat y}(\hat z)\leq 0$, then
 keeping $\tilde x$ fixed,  we then apply formula (\ref{eq4.12}) with $\hat
 y$, $\hat z$ in place of $y$, $z$, to the effect of yet further increase
 of the Poincar\'e coordinate of $\tilde z$ compared to that of $\hat
 z$ (and so of $z$ itself), measured within respective symmetric
 $y$-domains.
In case of $p_{\hat y,-\hat y}(\hat z)>0$ the image of point $z$
has already past the midpoint of the (varying) symmetric
$y$-domain interval while $y$-variable has been changing from $y$
to $\hat y$. Again, we then keep $\tilde x$ fixed, to move the
$y$-variable further, from $\hat y$ to $\tilde y$. This time,
application of formula (\ref{eq4.12})   can induce some decrease
in the Poincar\'e coordinate of the outcome -- the resulting point
$ p^{-1}_{\tilde y \tilde x}(c)$, with $c$ as in (\ref{pointc}),
can divide the  $y$-domain interval $(\tilde y,-\tilde y)$ in
smaller proportion than $\hat z$ did in $(\hat y,-\hat y)$.
Anyway, due to sharp inequality in (\ref{eq4.12}) for all $z$ such
that $p_{y,-y}(z)<0$, the midpoint could only be attained from the
other side. In other words, inequality (\ref{eq4.12}) guarantees
that the derivative of the induced $z$-movement, measured in the
respective Poincar\'e coordinates, is positive (and actually
bounded away from $0$) as long as the values assumed by the
$z$-variable are non-positive. Thus, in particular the value $0$
can be attained only once, and so if we put a point $\hat z$ with
$p_{y,-y}(\hat z)>0$ into the formula at the left-hand side of
inequality (\ref{eq4.12}), we necessarily  end up with a point on
the same side of $0$. Because the starting point $z$ was on the
other side, the lemma holds in this case too.
  \hfill$\Box$

  With lemma \ref{osemka} in place, we are in the position to state
  and prove the main result of this section.
\begin{thm}\label{thm4.1}
In the power-law family $f_a : x\mapsto -|x|^\alpha +a$, with
$\alpha >1$, there exists unique parameter $a=a(\alpha)$
corresponding to the kneading sequence $RLRRRLRC$.
 \end{thm}
{\it Proof}. In the course of the proof we make use of the tools
developed in section 2 of \cite{P}, where we pointed to some
consequences of homogeneity of the power-law mappings. In
particular, we had Lemma 2.1 there, asserting that for any two
points $q,\tilde q\in (0,1) $ one has
\begin{equation}\label{starylemat}
\overline h^-_{q,1}-\overline h^-_{\tilde q,1}=
 (p_{0,1}(h(q))-p_{0,1}(h(\tilde q)) )-(p_{0,1}(q)-p_{0,1}(\tilde q)
 ) .
\end{equation}
Speaking colloquially,  identity (\ref{starylemat}) tells, that
when we move the endpoint of an interval $(1,q)$  in $(1,0)$
towards the critical point, then an extra gain in the Poincar\'e
coordinate, coming from the successive action of the power-law
map, is just enough to make up for the loss (measured in
non-euclidean metric in $(1,q)$ and $(1,\tilde q)$ respectively)
suffered because of the simultaneously increased strength of the
limit non-euclidean push towards that moving endpoint. Other
propositions and lemmas of section 2 of \cite{P} served to
establish, that this limit situation, corresponding to Poincar\'e
coordinate close to $-\infty$, is essentially the worst possible,
and when we consider an interior point of a definite Poincar\'e
coordinate rather than the limit case, then the balance of gains
vs. losses is in our favor ("we are never in the red"). We will be
sending upon those properties when necessary, without reproducing
them in this paper.

Proceeding similarly to what we did in the proof of Theorem
\ref{thm3.1}, we split the procedure leading from $4_a$ to $8_a\,
,$ and respectively from $4_{a'}$ to $8_{a'} ,$ into several
steps. First, we increase $\bar a$ to $\bar a'$.  Theorem
\ref{thm3.1} yields, in particular,    that
$p_{0,-2_{a'}}(4_{a'})-p_{0,-2_a}(4_a)>0$. Next, we act upon
$4_{a'} $, and $ 4_a$, by the map $h$, and under the action of
$\overline h$ the above discrepancy gets enlarged. This is so,
because due to homogeneity,  we may for the purpose of performing
this step, tentatively  set each of the endpoints, $-2_{a'}$ and
respectively $-2_a$, at $1$. Then each of the Poincar\'e
coordinates $p_{0,-2_{a'}}(4_{a'})$, $\, p_{0,-2_a}(4_a)$, is
transformed by same, fixed negative Schwarzian map ${\overline
h}_{0,1}$. In the following step, we turn each of the intervals
$(0, h(2_{a'}))$, $\, (0, h( 2_a))\,$ over, and stretch them onto
$ (1, 3_{a'})\,$ and respectively $(1,3_a)$. The image of $4_{a'}$
is $\, 5_{a'}$, and by the so far described steps, it is clear
that
$p_{1,3_{a'}}(5_{a'})-p_{1,3-a}(5_a)>p_{0,-2_{a'}}(4_{a'})-p_{0,-2_a}(4_a)$.
By the truncation  argument from the proof of Theorem
\ref{thm3.1}, we know that $p_{1,0}(3_{a'})-p_{1,0}(3_a)>\bar
a-a\,$. In particular, the nonlinearity of the interval
$(1,3_{a'})$ is larger than that of $(1,3_a)$. Now, we act by the
homogeneous map $h$ again. Notice, that unlike in the case of
${\overline h}_{0,2_{a'}}$, this time the mapping ${\overline
h}_{1,3_{a'}}$ does not coincide with ${\overline h}_{1,3_a}$.
Anyway, we can still  claim that in this step the discrepancy of
the respective Poincar\'e coordinates grows again, i.e
\begin{equation}\label{kroczek}
 p_{1,h(3_{a'})}(h(5_{a'}))-p_{1,h(3_a)}(h(5_a))>p_{1,3_{a'}}(5_{a'})-
  p_{1,3_a}(5_a)\, .
\end{equation}
To this end, we invoke Propositions 2.5 and 2.4 of \cite{P}. From
the former, it follows that the strength of the non-euclidean push
generated by $h$ restricted to some domain,   is a monotone
function of the nonlinearity of that domain, when measured for a
fixed Poincar\'e coordinate within the varying domain. From the
latter, we derive that when the domain stays fixed, the strength
of the non-euclidean push of $h$ is monotone in the Poincar\'e
coordinate of the argument. We have noticed already that the
nonlinearity of $(1,3_{a})$ is increasing in parameter $\bar a$,
and also that $p_{1,3_{a'}}(5_{a'})>p_{1,3_a}(5_a)\, $, so the
principle of monotone behaviour of the strength of non-euclidean
push can be applied to the triples of points we consider. This
immediately implies the desired increase in the discrepancy of
appropriate Poincar\'e coordinates, as stated in (\ref{kroczek}).

 Making the next step, we turn the obtained triples
 $(1,h(5_{a'}),h(3_{a'}))$ and $\, (1,h(5_a),h(3_a))$ over, onto
 $(2_{a'},6_{a'},4_{a'}) $, and respectively onto $(2_a,6_a,4_a)$,
and then truncate them at the critical point $0$. In the proof of
Theorem \ref{thm3.1}, as well as in a step above, we were
satisfied to ascertain that this truncation  increases the
Poincar\'e coordinates discrepancy, which in current step would
yield $
p_{2_{a'},0}(6_{a'})-p_{2_a,0}(6_a)>p_{2_{a'},4_{a'}}(6_{a'})-p_{2_a,4_a}(6_a)
$, because by Theorem \ref{thm3.1} we know that
$p_{2_{a'},-2_{a'}}(4_{a'})>p_{2_a,-2_a}(4_a) $.
 To proceed further, one more observation is needed. It is fairly
clear that we have following lower bound on the increase of the
Poincar\'e coordinates discrepancy, generated by the cut-off at
$0$:
\begin{equation}\label{logarytm}
p_{2_{a'},0}(6_{a'})-p_{2_a,0}(6_a)>p_{2_{a'},4_{a'}}(6_{a'})-p_{2_a,4_a}(6_a)
+ \ln\frac{4_{a'}-2_{a'}}{4_a-2_a}\, .
\end{equation}
The equality in (\ref{logarytm}) is the limit case, attained for
infinitesimally short intervals placed at the left-hand endpoints,
i.e. when $p_{2_a,4_a}(6_a)\rightarrow -\infty$ and simultaneously
$p_{2_{a'},4_{a'}}(6_{a'})\rightarrow -\infty$. For
non-infinitesimal intervals  satisfying
$p_{2_{a'},4_{a'}}(6_{a'})>p_{2_a,4_a}(6_a)$, the same argument as
in the proof of Theorem \ref{thm3.1} obviously yields sharp
inequality in (\ref{logarytm}), and so the growth of the
discrepancy gained in the cut-off step is strictly larger than the
logarithmic term.

In the following step we once more  act by homogeneous map $h$,
and because ${\overline
 h}_{0,2_{a'}}$ coincides with ${\overline h}_{0,2_a}$, the same
 argument as before gives
\begin{equation}\label{jednrodne}
p_{h(2_{a'},0)}(h(6_{a'}))-p_{h(2_a),0}(h(6_a))>p_{2_{a'},0}(6_{a'})-
p_{2_a,0}(6_a)\, .
\end{equation}
This adds yet an extra amount to the discrepancy we consider. We
again turn the intervals $(h(2_{a'}),0) $ and $(h(2_a),0)$ over
and stretch them onto $(1,3_{a'})$ and $(1,3_a) $, with $6_{a'}$
going onto $7_{a'}$ and $6_a$ going onto $7_a$ respectively. It
remains to examine what happens in the last step, when we act by
the respective (non-coinciding!) restrictions of $h$ to the
obtained intervals, before we eventually return onto
$(2_{a'},4_{a'}) $ and onto $(2_a,4_a) $ by linear rescaling. This
is what we need Lemma \ref{osemka} for. In what follows we verify
its assumptions are fulfilled in our setting.

 In this last step we perform, the
strength of the non-euclidean push induced by $h{\left
|\right.}_{(1,3_{a'})}\,$, measured at $7_{a'}\,$, can be greater
than the respective strength of $h{\left |\right.}_{(1,3_a )}\,$
at $7_a$. This means that the discrepancy accumulated in all the
so far steps can now diminish. However, the identity
(\ref{starylemat}) provides a bound from the above on the amount
of possible loss. To see this, we recall that
\begin{equation}\label{inicjalnyprzyrost}
\bar a'-\bar a=p_{0,-1}(2_{a'})-p_{0,-1}(2_a)
<p_{1,0}(3_{a'})-p_{1,0}(3_a)\, ,
\end{equation}
and according to \ref{starylemat} we have
\begin{equation}\label{ostatnikrok}
(p_{2_{a'},1}(4_{a'})-p_{2_a,1}(4_a))-(p_{1,0}(3_{a'})-p_{1,0}(3_a))=
({\overline h }^-_{1,3_a}-{\overline h}^-_{1,3_{a'}} )
\end{equation}
We know that $p_{3_{a'},1}(7_{a'})>p_{3_a,1}(7_a)$ and the
interval $(1,3_{a'})$ has larger nonlinearity than $(1,3_a)$, so
we are in a position to invoke Propositions 2.5 and 2.4 of
\cite{P} once more. By them we have
\begin{equation}\label{push}
(p_{3_{a'},1}(7_{a'})-p_{3_a,1}(7_a) )-(p_{4_{a'},2_{a'}}(8_{a'})-
p_{4_a,2_a}(8_a)  )<({\overline h }^-_{1,3_a}-{\overline
h}^-_{1,3_{a'}} ) .
\end{equation}
The inequalities (\ref{logarytm}), (\ref{inicjalnyprzyrost}) and
(\ref{push}) put together, provide for fulfillment of condition
({\it iii}) of Lemma \ref{osemka}, with the points $x$, $y$ and
$z$ assuming values $2_a$, $4_a$ and $8_a$, as indicated before
the statement of the lemma. Now the claim of Theorem \ref{thm4.1}
follows directly from Lemma \ref{osemka}, and so we are done.
\hfill$\Box$

  We complete this section explicitly  recording one extra property, which
 we actually proved along the way.
 Denote the variable $\tau =p_{4_a, -4_a}(8_a)$
and let $\gamma = p_{2_a,-2_a}(4_a)$. From the proofs of Theorem
\ref{thm4.1} and Lemma \ref{osemka} there immediately follows
\begin{cor}
The function $\gamma =\gamma(\tau) : {\mathbb{R}}_-\to{\mathbb{R}
}_+ \, $ is strictly increasing in $\tau$, with the derivative
$\gamma'(\tau)$ bounded away from $0$ and $+\infty$.
\end{cor}
\hfill$\Box$
\section{$RLRLRLRL\ldots$}\label{rlrl}
In this section we let the parameter $\bar a$ vary in  a range
such that the kneading sequence is $RLRL\ldots\,$. From Theorem
\ref{thm3.1} it follows immediately that the range of admissible
$\bar a$'s $\,$ is  always a half-line $\,( -\infty, {\bar a}_1)$,
with the specific value of ${\bar a}_1$ depending on the critical
order $\alpha$.
 Upholding
the normalization $0_a=0 \,$, $1_a=1$ we have set before, this
means the post-critical orbit begins with $2_a\in (0,-1)\, $,
$3_a\in (1,0) $ and $4_a\in (2_a,0) $. Then, we get two sequences
of nested intervals, the odds: $(1,3_a)$, $(3_a,5_a)$,
$(5_a,7_a)\ldots\,$, and the evens: $(0,2_a)$, $(2_a,4_a)$,
$(4_a,6_a)\ldots\,$.
 In
terms of multipliers, we either have a period $2$ orbit with
negative multiplier, or  this periodic orbit had turned into a
repeller  and, by bifurcation, there was born a period $4$
attracting periodic orbit with positive multiplier. In what
follows, we shall  see that the ratios of consecutive intervals
within each of the two decreasing families are functions strictly
monotone in parameter $\bar a$. Moreover, the initial increase of
the parameter, i.e. $\bar a'-\bar a$, does not eventually vanish,
but a definite part of it is preserved through all the steps. This
will further provide, with some extra work, for monotonicity of
the multipliers, also in the case of repelling period $2$ orbit.
This is a work in preparation. The remaining part of this paper is
devoted to the proof of the following claim.
\begin{thm}\label{nestedfamily}
For $\bar a\in (-\infty,{\bar a}_1)$ and for all non-negative
integers $n$, the ratio functions
\begin{equation}
r^n_e = \frac{|(2n+4)_a - (2n+2)_a|}{|(2n+2)_a - (2n)_a|} \ \ \ \
{\rm and} \ \ \ \ r^n_o = \frac{|(2n+5)_a - (2n +3)_a|}{|(2n+3)_a
- (2n+1)_a|}
  \end{equation}
are strictly increasing in $\bar a$.

Moreover, when   the parameter increases from $\bar a$ to $\bar
a'$, then for every  $n\in \mathbb{Z}_+$ the induced discrepancy
of the Poincar\'e coordinates satisfies
\begin{equation}
 p_{(n+2)_{a'},n_{a'}}((n+4)_{a'}) - p_{(n+2)_{a},n_{a}}((n+4)_{a})
 > (p_{1,-1}(2_a') - p_{1,-1}(2_a))\, .
\end{equation}
\end{thm}
{\it Proof}. As before, we divide the procedure into steps. Once
we cover the most delicate step, which  turns out to be the
passage from $(5_{a'},7_{a'})\,$ to $\, (6_{a'},8_{a'})$, we will
be in a position to continue inductively. We begin by moving the
initial point $2_a$ to a new position $2_{a'}$, with $\bar a'>\bar
a$. Then, by the truncation argument from the proof of Theorem
\ref{thm3.1}, we have
\begin{equation}\label{napoczatek}
 (p_{1,0}(3_{a'})-p_{1,0}(3_a))>\bar a' -\bar a=\Delta\bar a
>\Delta t  ,
\end{equation}
where we denoted $ \Delta t= (p_{1,-1}(2_a') - p_{1,-1}(2_a))$.
  Since we then act by
the homogeneous map $h$, by  (\ref{starylemat}) we get
\begin{equation}
p_{2_{a'},1}(4_{a'})-p_{2_a,1}(4_a)=(p_{1,0}(3_{a'})-p_{1,0}(3_a))+
( {\overline h}^-_{3_a,1} - {\overline h}^-_{3_{a'},1}   )\, .
\end{equation}
Passing from $3_{a'}$ to $4_{a'}$, we cannot directly apply the
truncation argument again, because in this step the Poincar\'e
coordinate $p_{1,2_{a'}}(0)$ of the cut-off point decreases (cf.
Proposition \ref{prop2.1}). That can be fixed by decomposing the
step in two, and simultaneous use of  Proposition \ref{prop2.4},
identity (\ref{starylemat}) and truncation. According to
(\ref{napoczatek}) and Proposition \ref{prop2.4},
\begin{equation}
p_{2_{a'},-2_{a'}}\left(p^{-1}_{2_{a'},1}(p_{1,0}(3_{a'}))
\right)- p_{2_{a},-2_{a}}\left(p^{-1}_{2_{a},1}(p_{1,0}(3_{a}) )
\right) > \Delta t .
\end{equation}
Truncation at $0$ obviously gives
\begin{equation}
p_{2_{a'},0}\left(p^{-1}_{2_{a'},1}(p_{1,0}(3_{a'})) \right)-
p_{2_{a},0}\left(p^{-1}_{2_{a},1}(p_{1,0}(3_{a}) ) \right) >
\Delta t .
\end{equation}
Then, to the Poincar\'e coordinate of the point
$(p^{-1}_{2_{a'},1}(p_{1,0}(3_{a'}))$, read in the domain
$(2_{a'},1) $, we add the extra gain of $( {\overline h}^-_{3_a,1}
- {\overline h}^-_{3_{a'},1} )$. The non-euclidean length of this
same extra interval, read in the domain $(2_{a'},0)$ rather than
in $ (2_{a'},1) $, is of course larger, because of truncation.
Thus
\begin{equation}\label{nanapoczatek}
p_{2_{a'},0}(4_{a'})-p_{2_a,0}(4_a)>\Delta t +( {\overline
h}^-_{3_a,1} - {\overline h}^-_{3_{a'},1}   )\, .
\end{equation}
Doing  the homogeneous mapping again, by (\ref{starylemat}) and
(\ref{nanapoczatek}) we get
\begin{eqnarray}\label{nananapoczatek}
p_{3_{a'},1}(5_{a'})-p_{3_{a},1}(5_{a})= {\overline
h}_{2_{a'},0}(p_{2_{a'},0}(4_{a'}))-{\overline
h}_{2_{a},0}(p_{2_{a},0}(4_{a}))>\nonumber
\\
>\Delta t +( {\overline h}^-_{3_a,1} - {\overline h}^-_{3_{a'},1})
+( {\overline h}^-_{4_a,2_a} - {\overline h}^-_{4_{a'},2_{a'}})\,
.
\end{eqnarray}
Now, similarly to the final step in the proof of Theorem
\ref{thm4.1}, we can argue that the so far acquired gain in the
Poincar\'e coordinate is enough to make up for possible losses in
the next two steps. This is fairly clear. The interval
$(3_{a'},1)$ has larger nonlinearity than $(3_a, 1)$, and
$p_{3_{a'},1}(5_{a'})>p_{3_a,1}(5_a)$, so Propositions 2.5 and 2.4
of \cite{P} do apply when we act by
${h\left|\right.}_{(1,3_{a'})}$ and ${h\left|\right.}_{(1,3_{a})}
.$ Therefore, in this step the discrepancy
$(p_{3_{a'},1}(5_{a'})-p_{3_a,1}(5_a))$ can only be diminished by
an amount smaller than $  ( {\overline h}^-_{3_a,1} - {\overline
h}^-_{3_{a'},1})$, yielding
\begin{equation}\label{nanananapoczatek}
p_{4_{a'},2_{a'}}( 6_{a'})-p_{4_{a},2_{a}}( 6_{a})>\Delta t + (
{\overline h}^-_{4_a,2_a} - {\overline h}^-_{4_{a'},2_{a'}}),
\end{equation}
By (\ref{nanapoczatek}), the nonlinearity of $(4_{a'}, 2_{a'})$
 is
larger than that of $(4_{a}, 2_{a})$, and also $p_{4_{a'},2_{a'}}(
6_{a'})>p_{4_{a},2_{a}}( 6_{a}) $. Thus, when we act by
 ${h\left|\right.}_{(4_{a'},2_{a'})}$, and respectively by
 ${h\left|\right.}_{(4_a,2_{a})}$, we certainly do not lose more
 than $ (
{\overline h}^-_{4_a,2_a} - {\overline h}^-_{4_{a'},2_{a'}}) $ in
the outgoing discrepancy.
 Hence, by ( \ref{nanananapoczatek})
 \begin{equation}\label{pip}
p_{5_{a'},3_{a'}}(7_{a'})-p_{5_a,3_a}(7_a)>\Delta t .
 \end{equation}
We can now make a shortcut towards completion of the current
cycle. The nonlinearity of $ (3_{a'},1 )$ is larger than that of $
( 3_a,1) $ and $ p_{3_{a'},1}(5_{a'})>p_{3_a,1}(5_a) $, which in
turn gives that the nonlinearity of $ ( 5_{a'},3_{a'}) $ is larger
than that of $(5_a,3_a ) $. Also $
p_{5_{a'},3_{a'}}(7_{a'})>p_{5_a,3_a}(7_a) $, so we can apply the
argument about monotonicity of the strength of the non-euclidean
push, which we recalled in the proof of Theorem \ref{thm4.1},
immediately arriving at
\begin{equation}\label{pippip}
p_{6_{a'},4_{a'}}(8_{a'})-p_{6_a,4_a}(8_a)>p_{5_{a'},3_{a'}}(7_{a'})-p_{5_a,3_a}(7_a)>\Delta
t .
\end{equation}
However, the above argument alone turns out to be insufficient,
when we want to do further iterates. To obtain an inequality which
we could use inductively at all steps, we need more subtle
understanding at this particular stage of our procedure. Here we
go.

 From (\ref{napoczatek}) and Proposition \ref{prop2.5} it follows
that
\begin{equation}
p_{p^{-1}_{3_{a'},1}(p_{3_a,1}(5_a)),0}(3_{a'})-
 p_{5_a,0}(3_a)>\Delta\bar a\, ,
\end{equation}
so by  $\, p_{3_{a'},1}(5_{a'})>p_{3_{a},1}(5_{a})$ we have
\begin{equation}
p_{5_{a'},0}(3_{a'})-p_{5_a,0}(3_a)>\Delta\bar a .
\end{equation}
By the same argument applied to $(5_{a'},3_{a'})$ rather than
$(1,3_{a'})$, we get
 \begin{equation}\label{ku}
p_{ p^{-1}_{5_{a'},3_{a'}}( p_{5_a,3_a}(7_a) ),0
}(3_{a'})-p_{7_a,0}(3_a)>\Delta\bar a .
\end{equation}
Now we do the homogeneous mapping $h$, and rescale the image onto
$(1,2_{a'})$. The image of $3_{a'}$ is $4_{a'}$, and by a version
of the   truncation argument alike that used before in the step
leading from $(1,3_{a'})$ to $(2_{a'},4_{a'})$, we use
\begin{equation}
p_{2_{a'},0}( p^{-1}_{2_{a'},1}(p_{1,0}(3_{a'})) )-p_{2_a,0}(
p^{-1}_{2_a,1}(p_{1,0}(3_a)) )>\Delta t
\end{equation}
 and (\ref{ku}) to get
  \begin{equation}\label{kuku}
 p_{ p^{-1}_{2_{a'}1}( {\overline h}_{1,0}( p_{1,0}( p^{-1}_{5_{a'},3_{a'} }( p_{5_a,3_a}(7_a ) ) ) )
 ),0}(4_{a'})-p_{8_a,0}(4_a)>\Delta t .
 \end{equation}
This is so, because (\ref{ku}) implies
\begin{equation}\label{fiu}
p_{h( p^{-1}_{5_{a'},3_{a'}}( p_{5_a,3_a}(7_a) )),0
}(h(3_{a'}))-p_{h(7_a),0}(h(3_a))>\Delta\bar a ,
\end{equation}
 and when we consider the interval $(c,d)$, where $d=p^{-1}_{2_{a'}1}( p_{2_a,1}(4_a)+\Delta\bar
 a)$,
and the point $c$ is defined so that
 \begin{equation}\label{fiufiu}
p_{c,1}(d)-p_{8_a,1}(4_a)=\Delta\bar a
 \end{equation}
 then, according to
Proposition \ref{prop2.5} applied to the domain $(1,2_{a'})$ in
place of $(0,1)$, and with the points $d$ and $c$ singled out, we
see that point $c $ divides the interval $(d,2_{a'})$ at the same
proportion as $8_a$ divided $(4_a,2_a)$. Re-applying Proposition
\ref{prop2.5} to the domain $(0,2_{a'})$ with the same singled out
pair of points, we further see that
\begin{equation}\label{fiut}
p_{c,0}(d)>p_{8_a,0}(4_a)+\Delta t ,
\end{equation}
 because by Proposition
\ref{prop2.4} $\, p_{2_{a'},0}(d)-p_{2_a,0}(4_a)>\Delta t$.
Recalling (\ref{fiu}) and taking into account that
$p_{2_{a'},1}(4_{a'})>p_{2_{a'},1}(d) $, which in turn gives
$p_{1,4_{a'}}(0)>p_{1,d}(0)$, we can now do the standard
truncation argument, cutting-off at $0$ to arrive at (\ref{kuku}).

 This formula could do for the iterative procedure if we cared
 only for some, indefinite growth. To obtain definite growth, claimed
 in the statement of Theorem \ref{nestedfamily}, we need to work harder.

 In the next step of the proof, we will  see that the extra amount
of $\Delta t$ in formula (\ref{kuku}) allows us to move $7_a$
towards the endpoint by at least that much. To this end, we again
consider the interval $(5_{a'},3_{a'})$, but this time the point
within we single out, is point $e$ determined by
\begin{equation}\label{punkte}
p_{5_{a'},3_{a'}}(e)=p_{5_a,3_{a}}(7_a)+\Delta t.
\end{equation}
 From
(\ref{ku}), using Proposition \ref{prop2.1} with points $0$,
$3_{a'}$ and $e$ in place of $1$, $0$ and $x$ respectively, or by
a direct check, one gets
\begin{equation}\label{tf}
p_{e,0}(3_{a'})-p_{7_a,0}(3_a)>\Delta \bar a -\Delta t .
\end{equation}
Doing the homogeneous mapping, we have
\begin{equation}\label{tfu}
p_{h(e),0}(h(3_{a'}))-p_{h(7_a).\,0}(h(3_a))>\Delta \bar a -\Delta
t .
\end{equation}
Again, we consider an interval $(f,d)$, where  $d$ has same
meaning as above, and point $f$ is defined by
\begin{equation}\label{fiufiufiu}
p_{f,1}(d)-p_{8_a,1}(4_a)=\Delta\bar a -\Delta t .
 \end{equation}
From (\ref{fiufiu}) and (\ref{fiufiufiu}), it follows by
Proposition \ref{prop2.1} that
$p_{d,-\infty}(c)-p_{d,-\infty}(f)=\Delta t$, and again by this
same proposition $p_{c,0}(d)-p_{f,0}(d)=\Delta t$. Hence, by
(\ref{fiut}), we have $p_{f,0}(d)>p_{8_a,0}(4_a)$. This, and
(\ref{tfu}) lead to
\begin{equation}\label{tfutfu}
p_{p^{-1}_{6_{a'},4_{a'} }({\overline h}_{5_{a'},3_{a'}
}(p_{5_a,3_a}(7_a)+\Delta t ) ),0 }(4_{a'})>p_{8_a,0}(4_a).
 \end{equation}

  We can describe what we have found so far in the following way.
  We move the parameter up, from $\bar a$ to $\bar a'$. In the odd
  family, we see $3_a$ moving to $3_a'$ by more than $\Delta\bar
  a$. Consequently,  the non-euclidean coordinate of $5_a$ vary,   within
  its dynamically determined base interval,  by at least
  $\Delta t$, plus an additional increment which is sufficient to
  make up for the increased -- due to larger nonlinearity of the new
  new domain intervals -- strength of the non-euclidean push
  backwards. In the next odd return we  do not let $7_a$ move
  all the way to its new position $7_{a'}$ at once. Instead, we
  first only add $\Delta t$ to its Poincar\'e coordinate. This
  corresponds to starting from the point $e$ in the already fully
  enlarged domain $(5_{a'},3_{a'})$, rather than from $7_{a'}$. We
  have just seen that not only is the nonlinearity of $(e,3_{a'})$
  larger than that of $(7_a, 3_a)$, but the nonlinearity of
  $(\hat e,4_{a'})$ is larger than that of $(8_a,4_a)$ also. Here $\hat e$
  is the dynamical successor of $e$ on the even side. This latter
  estimate from the below on the the nonlinearity, turns out to be
  fundamental for the prospective iterates.

Recall we defined $e$ by (\ref{punkte}) so as to have $
p_{5_{a'},3_{a'}}(e)=p_{5_a,3_{a}}(7_a)+\Delta t$. The same way we
derived (\ref{pippip}) from (\ref{pip}) we also get
  \begin{equation}\label{bb}
p_{6_{a'},4_{a'}}(\hat e)-p_{6_a,4_a}(8_a)>\Delta t .
\end{equation}
 This will be
needed, when it comes to definite growth in both odd and even
family. But now, for points $e$ and $\hat e$ we have stronger
input: in both cases, we know that the nonlinearity of the
remaining part of the  base interval increased. Therefore, we will
now be able to proceed pretty much like in the initial step, that
led from $(1,3_{a'})$ to $(2_{a'},4_{a'})$, rather than use the
earlier described shortcut.  Similarly to that initial step, we
again want to know that the surplus exceeding $\Delta t$  in
(\ref{bb}) will make up for possible loss, inflicted by increased
nonlinearity of $( 5_{a'},7_{a'})$, upon next return to $
(5_{a'},7_{a'})$. However, we have to overcome a serious obstacle.
Formula (\ref{starylemat}) we previously used to that goal, holds
true only so long as the critical point is the endpoint
corresponding to non-euclidean $+\infty$. This is of course not
the case for $(5_a,3_a)$, nor for all other intervals in our odd
and even families, except for the initial ones. For intervals not
bounded by the critical point, we only know monotonicity of the
strength of non-eucliean push and this, in general, does not give
control over an amount of the gain in Poincar\'e coordinates
discrepancy. Composing ${\overline h}$ mappings over two
arbitrary, successive domains, yet worsens the the situation.
Fortunately, all this can be fixed with (\ref{tf}) and
(\ref{tfutfu}) in place. Increased nonlinearity of that  part of a
domain interval  which bounds us away from the endpoint, provides
an effective
 replacement for the critical endpoint. In particular, we will see
 that the gain in the non-euclidean coordinates discrepancy is
 even better than that in
 formula  (\ref{starylemat}) . This is why we have striven for
those nonlinearity inequalities. As soon as we are over with the
part which takes $7_a$ to $e$, the remaining part, in which we
move $e$ further to $7_{a'}$, will require only an easy estimate.
All the above holds true for even successors, $\hat e$ and
eventually $8_{a'}\,$, as well. With one extra observation to
make, we will be able to do arbitrarily long iterates, preserving
the $\Delta t$ discrepancy all along the way.

 To carry out the  above described  strategy, we recall that in
 Proposition 2.2 of \cite{P} we gave an explicite formula for the
 strength of non-euclidean push, which  turns out to be
 \begin{equation}
 |p_{r,s}(\varphi
 (x))-p_{p,q}(x)|=|{\overline\varphi}^-_{x,q}+{\overline\varphi}^+_{p,x}|
 \end{equation}
We also noticed there, that for the homogeneous mapping $h$
restricted to some interval, the quantities ${\overline h}^-$ and
${\overline h}^+$ depend solely on the nonlinearity of that domain
interval. By monotonicity of the strength of the non-euclidean
push as a function of the nonlineatity of the domain, also the
limit values, ${\overline h}^-$ and ${\overline h}^+$,  behave
monotonically.  By all the above, taking (\ref{tf}) into account,
we have
\begin{eqnarray}\label{pa}
p_{6_{a'},4_{a'}}(\hat e)-p_{6_a,4_a}(8_a)=
  {\overline h}_{5_{a'},3_{a'}
}(p_{5_a,3_a}(7_a)+\Delta t )  -p_{6_a,4_a}(8_a)= \nonumber\\
 p_{h(5_{a'}),h(3_{a'}) }(h(e))-p_{h(5_a),h(3_a)}(h(7_a)) >\Delta t +({\overline
 h}^+_{5_{a'},e}-{\overline h}^+_{5_{a},7_a})
\end{eqnarray}
The sign at the superscript of ${\overline h}$ in (\ref{pa})
depends only on an orientation of the domain, so (\ref{pa})
provides  a better estimate than we could derive from
(\ref{starylemat}), if the endpoint $3_{a'}$ coincided with the
critical point. Doing the successive $h$-map step on the even
side, because of (\ref{tfutfu}), we get in the same way
\begin{equation}\label{prr}
p_{h(6_{a'}),h(4_{a'}) }(h(\hat e))-p_{h(6_a),h(4_a)}(h(8_a))
>\Delta t +({\overline
 h}^+_{5_{a'},e}-{\overline h}^+_{5_{a},7_a})+
 ({\overline h}^+_{6_{a'},\hat e} - {\overline h}^+_{6_a,8_a} ) .
\end{equation}
These are formulas analogous to (\ref{nanapoczatek}) and
({\ref{nananapoczatek}}), and what we  want now, is a similar
estimate where the input is $7_{a'}$ and $8_{a'}$, rather than $e$
and $\hat e$. To move from $e$ to $7_{a'}$ we could simply invoke
Lemma 2.4 of \cite{P}. However,  there is no generalization  of
that lemma which could be used over two unrelated domains. We need
to be a bit more careful, and use the dynamical relation between
an interval and its image.  Doing the mapping $h$, by homogeneity
and (\ref{starylemat}) we have
 \begin{equation}\label{pr}
p_{h(5_{a'}),0}(h(7_{a'}))-p_{h(5_{a'}),0}(h(e))=(p_{5_{a'},0}(7_{a'})-p_{5_{a'},0}(e))+
({\overline h}^+_{5_{a'},7_{a'}} -{\overline h}^+_{5_{a'},e}) .
\end{equation}
Before we do another mapping $h$, we take the image over onto
$(1,-\infty)$, so that $h(3_{a'})$ goes onto $(4_{a'})$, and cut
off at $0$. Because of this truncation
\begin{equation}
p_{6_{a'},0}(8_{a'})-p_{6_{a'},0}(\hat
e)>(p_{5_{a'},0}(7_{a'})-p_{5_{a'},0}(e))+ ({\overline
h}^+_{5_{a'},7_{a'}} -{\overline h}^+_{5_{a'},e}) .
\end{equation}
 Now, acting
by  homogeneous map, we get
\begin{eqnarray}
p_{h(6_{a'}),0}(h(8_{a'}))-p_{h(6_{a'}),0}(h(\hat
e))>(p_{5_{a'},0}(7_{a'})-p_{5_{a'},0}(e))+\nonumber\\
+({\overline h}^+_{5_{a'},7_{a'}} -{\overline h}^+_{5_{a'},e})
 + ({\overline
h}^+_{6_{a'},8_{a'}} -{\overline h}^+_{6_{a'},\hat e})
\end{eqnarray}
 We neglect a positive summand $ (p_{5_{a'},0}(7_{a'})-p_{5_{a'},0}(e)) $
and truncate at $ h(4_{a'}) $ to arrive at
\begin{equation}\label{prrprr}
p_{h(6_{a'}),h(4_{a'})}(h(8_{a'}))-p_{h(6_{a'}),h(4_{a'})}(h(\hat
e))>({\overline h}^+_{5_{a'},7_{a'}} -{\overline h}^+_{5_{a'},e})
 + ({\overline
h}^+_{6_{a'},8_{a'}} -{\overline h}^+_{6_{a'},\hat e})
\end{equation}
Similarly, neglecting a positive summand at (\ref{pr}),  followed
by cutting off at $h(3_{a'})$ leads to
\begin{equation}\label{prpr}
p_{h(5_{a'}),h(3_{a'})}(h(7_{a'}))-p_{h(5_{a'}),h(3_{a'})}(h(e))>
({\overline h}^+_{5_{a'},7_{a'}} -{\overline h}^+_{5_{a'},e}) .
\end{equation}
The inequalities (\ref{pa}) and (\ref{prpr}), in conjunction with
(\ref{prr}) and (\ref{prrprr}), give
\begin{equation}\label{pu}
 p_{h(5_{a'}),h(3_{a'}) }(h(7_{a'}))-p_{h(5_a),h(3_a)}(h(7_a)) >\Delta t +({\overline
 h}^+_{5_{a'},7_{a'}}-{\overline h}^+_{5_{a},7_a})
\end{equation}
and also
\begin{equation}\label{prrprrprr}
p_{h(6_{a'}),h(4_{a'}) }(h(8_{a'}))-p_{h(6_a),h(4_a)}(h(8_a))
>\Delta t +({\overline
 h}^+_{5_{a'},7_{a'}}-{\overline h}^+_{5_{a},7_a})+
 ({\overline h}^+_{6_{a'},8_{a'}} - {\overline h}^+_{6_a,8_a} ) ,
\end{equation}
which are the desired estimates. Now, by the same argument which
led from (\ref{nanapoczatek}) and (\ref{nananapoczatek}), through
(\ref{nanananapoczatek}) to (\ref{pip}), we can see that
(\ref{pu}) and (\ref{prrprrprr}) imply
\begin{equation}\label{qu}
p_{9_{a'},7_{a'}}(11_{a'})-p_{9_a,7_a}(11_a)>\Delta t .
\end{equation}
In the same way we obtained (\ref{pippip}) from (\ref{pip}), we
can also derive
\begin{equation}\label{ququ}
p_{10_{a'},8_{a'}}(12_{a'})-p_{10_a,8_a}(12_a)>\Delta t .
\end{equation}
We have completed the second cycle. Those  were necessary to
initialize the inductive procedure. We are now in a position to do
the final argument, which can be used repeatedly. We believe that,
because all the elaborate notation of the first two cycles is
already in place, it will be more instructive to present this
argument in detail as the next cycle, rather than in general
terms. It will be obvious that what we do, is tantamount to the
inductive step.

 We pick a point $\epsilon\in (9_{a'},7_{a'})$,
such that
\begin{equation}
p_{9_{a'},7_{a'}}(\epsilon)=p_{9_a,7_a}(11_a)+\Delta t .
\end{equation}
We will prove that the nonlinearity of $(\epsilon ,7_{a'})$ is
larger than that of $(11_a,7_a)$, and simultaneously  the
nonlinearity of $(\hat\epsilon ,8_{a'}) $ is larger than that of
$(12_a, 8_a)$. As before, $\hat\epsilon$ stands for the dynamical
successor of $\epsilon$ on the even side. This will permit to
bypass the non-critical endpoint obstacle in the next cycle, the
way we did earlier, with $e$ and $\hat e$. To show this
nonlinearity increase, we proceed in several steps. First, in
$(5_{a'},e)$ we find a point $\beta
=p^{-1}_{5_{a'},e}(p_{5_a,7_a}(9_a))$. Then, in $(\beta ,e)$ we
find $\varepsilon$, such that $p_{\beta
,e}(\varepsilon)=p_{9_a,7_a}(11_a)+\Delta t$. We use the the fact
that Poincar\'e coordinate of the point $e$, compared to that of
$7_a$, is already moved by $\Delta t$ towards the endpoint, to
ascertain that the nonlinearity of $(\varepsilon ,e)$ is larger
than nonlinearity of $(11_a,7_a)$. This is so, because according
to Proposition \ref{prop2.5}, $p_{\beta
,3_{a'}}(e)=p_{9_a,3_a}(7_a)+\Delta t$, and the nonlinearity of
$(e,3_{a'})$ is larger than nonlinearity of $(7_a,3_a)$; we have
\begin{equation}
\frac{|(\beta ,e)|}{|(e,3_{a'})|}=(\exp{\Delta
t})\frac{|(9_a,7_a)|}{|(7_a,3_a)|}>(\exp{\Delta t})\frac{|(\beta,
\delta )|}{|(\delta ,e)|} ,
\end{equation}
where $\delta\in (\beta ,e)$ is a point such that $\frac{|(\delta
,e )|}{|(e,3_{a'})|}=\frac{|(11_a,7_a)|}{|(7_a,3_a)|}$. Thus,
  $p_{e,\beta }(\varepsilon)>p_{e,\beta}(\delta)$, and
  consequently $p_{\varepsilon ,3_{a'}}(e)>p_{\delta ,3_{a'}}(e)\,$
  or, in other words, $\frac{|(\varepsilon, e)|}{|(e,3_{a'})|}>\frac{|(\delta ,e)|}
  {|(e,3_{a'})|}$. Since $(e,3_{a'})$ has larger nonlinearity than
  $(7_a,3_a)$, the nonlinearity of $(\varepsilon ,e)$ must be larger
than nonlinearity of $(11_a,7_a)$. Next, we do the mapping $h$ and
 consider the situation on the even side.  The interval
 $(\varepsilon,3_{a'})$ has larger nonlinearity than $(11_a,3_a)$
 and $p_{\varepsilon ,3_{a'}}(e)>p_{11_a,3_a}(7_a)$, so by
 principles  of monotonicity of the strength of non-euclidean push
 in nonlinearty of the domain, as well as in the coordinate of the
 point, the action of ${\overline h}_{\varepsilon ,3_{a'}}$ makes
$$p_{\hat\varepsilon ,4_{a'}}(\hat e)>p_{12_a,4_a}(8_a) ,$$
where $\hat\varepsilon$ is the dynamic successor of $\varepsilon$.
Because $(\hat e, 4_{a'})$ has larger nonlinearity than
$(8_a,4_a)$, it follows that the nonlinearity of $(\hat\varepsilon
,\hat e )$ is also larger than that of $(12_a,8_a)$.
 By the same two principles applied to $p_{\beta
 ,e}(\varepsilon)$,
 we get $p_{\hat \beta ,\hat e}(\hat\varepsilon)>p_{10_a,8_a}(12_a)+\Delta
 t$, but because of the proved nonlinearity increases, we can also
 claim that
 $$p_{h(\beta),h(e)}(h(\varepsilon))-p_{h(9_a),h(7_a)}(h(11_a))>\Delta t+(
 {\overline h}^+_{\beta ,\varepsilon}-{\overline h}^+_{9_a,11_a}  )\,\ \ \ \ \ \ \ \
 {\rm and} $$
 $$ p_{h(\hat\beta),h(\hat e)}(h(\hat\varepsilon))-p_{h(10_a),h(8_a)}(h(12_a))>\Delta t+(
 {\overline h}^+_{\beta ,\varepsilon}-{\overline h}^+_{9_a,11_a}
 )+(  {\overline h}^+_{\hat\beta ,\hat\varepsilon}-{\overline h}^+_{10_a,12_a}
 ) .
 $$
 We are through
with the first part of the inductive step. Now, our immediate plan
is to move $e$ to $7_{a'}$, then $\beta$ up to $9_{a'}$, and to
replace $\varepsilon$ by $\epsilon$, keeping all the  above gains
untouched,  both on the odd and on the even side. Having done all
that, we will easily be able to move $\epsilon$ to $11_{a'}$, to
complete the procedure.

  Denote by $\lambda$ the point determined by
  $p_{5_{a'},7_{a'}}(\lambda)=p_{5_a,7_a}(9_a)$, and let
 $\tau\in (\lambda ,7_{a'})$ be such, that $p_{\lambda ,7_{a'}}(\tau)=
 p_{9_a,7_a}(11_a)+\Delta t$. The interval $(5_{a'},7_{a'})$ has
 larger nonlinearity than $(5_{a'},e)$, so consequently $(\lambda
 ,7_{a'})$ has larger nonlinearity than $(\beta ,e)$, and $(\tau ,7_{a'})$
 has larger nonlinearity than $(\varepsilon ,e)$. Thus $(h(\tau)
 ,h(7_{a'}))$ has larger nonlinearity than $(h(\varepsilon) , h(e))$.
The distance of $8_{a'}$  to the critical point $0$ is smaller
than similar distance for  the point $\hat e$, so   the truncation
argument after cutting of at $0$, implies that $(\hat\tau
,8_{a'})$ has larger nonlinearity than $(\hat\varepsilon ,\hat e)$
and, in turn, larger than $(12_a,8_a)$. Again, we increase the
intervals in question, choosing $9_{a'}$ in place of $\lambda$,
and replacing $\tau$ by $\epsilon$. Then, of course, $(\epsilon
,7_{a'})$ has yet larger nonlinearity, so
$(h(\epsilon),h(7_{a'}))$ has larger nonlinearity than $(h(\tau)
 ,h(7_{a'}))$ and , after truncation, $(\hat\epsilon , 8_{a'})$ has
 larger nonlinearity than $ (\hat\tau
,8_{a'}) $.  It immediately implies
\begin{equation}
p_{h(9_{a'}),h(7_{a'})}(h(\epsilon))-p_{h(9_a),h(7_a)}(h(11_a))>\Delta
t+(
 {\overline h}^+_{9_{a'} ,\epsilon}-{\overline h}^+_{9_a,11_a}  )
 ,
\end{equation}
and
\begin{equation}
 p_{h(10_{a'}),h(8_{a'})}(h(\hat\epsilon))-p_{h(10_a),h(8_a)}(h(12_a))>\Delta t+(
 {\overline h}^+_{9_{a'} ,\epsilon}-{\overline h}^+_{9_a,11_a}
 )+(  {\overline h}^+_{\hat\beta ,\hat\varepsilon}-{\overline h}^+_{10_a,12_a}
 ) .
 \end{equation}
  This is what we aimed at.
 By the same argument that  earlier let us  replace  $e$ by
 $7_{a'}$ and $\hat e$ by $8_{a'}$ to derive formulas (\ref{pu}) and (\ref{prrprrprr}), we
  can now replace $\epsilon$ by $11_{a'}$ and ${\hat\epsilon}$ by
  $12_{a'}$, arriving at
\begin{equation}\label{tup}
 p_{h(9_{a'}),h(7_{a'}) }(h(11_{a'}))-p_{h(9_a),h(7_a)}(h(11_a)) >\Delta t +({\overline
 h}^+_{9_{a'},11_{a'}}-{\overline h}^+_{9_{a},11_a}) ,
\end{equation}
and
\begin{eqnarray}\label{tuptup}
p_{h(10_{a'}),h(8_{a'}) }(h(12_{a'}))-p_{h(10_a),h(8_a)}(h(12_a))
> \nonumber\\
 \Delta t +({\overline
 h}^+_{9_{a'},11_{a'}}-{\overline h}^+_{9_{a},11_a})+
 ({\overline h}^+_{10_{a'},12_{a'}} - {\overline h}^+_{10_a,12_a} ) .
\end{eqnarray}
Similarly to (\ref{qu}) and (\ref{ququ}), we also get
$$
p_{13_{a'},11_{a'}}(15_{a'})-p_{13_a,11_a}(15_a)>\Delta t ,\ \ \ \
{\rm and }\ \ \ \
p_{14_{a'},12_{a'}}(16_{a'})-p_{14_a,12_a}(16
_a)>\Delta t .
$$ This completes the inductive step. The claim of the theorem
follows immediately. \hfill $\Box$

\end{document}